\documentclass[12pt,a4paper]{article}
\usepackage{a4,amsmath,amssymb,amsfonts,amsthm}
\usepackage{amsrefs}
\usepackage[all]{xy}
\usepackage[small,nohug,heads=vee]{diagrams}
\diagramstyle[labelstyle=\scriptstyle]

\begin{document}
\bibliographystyle{plain}
 

\def\mR{\M{R}}           
\def\mZ{\M{Z}}           
\def\mN{\M{N}}           
\def\mQ{\M{Q}}       
\def\mC{\M{C}}  
\def\mG{\M{G}}



\def\deg{{\rm deg}}
\def\Spec{{\rm Spec}}
\def\rg{{\rm rg}}
\def\Hom{{\rm Hom}}
\def\Aut{{\rm Aut}}
 \def\Tr{{\rm Tr}}
 \def\Exp{{\rm Exp}}
 \def\Gal{{\rm Gal}}
 \def\End{{\rm End}}
 \def\det{{{\rm det}}}
 \def\Td{{\rm Td}}
 \def\ch{{\rm ch}}
 \def\che{{\rm ch}_{\rm eq}}
  \def\Spec{{\rm Spec}}
\def\Id{{\rm Id}}
\def\Zar{{\rm Zar}}
\def\Supp{{\rm Supp}}
\def\eq{{\rm eq}}
\def\Ann{{\rm Ann}}
\def\LT{{\rm LT}}
\def\Pic{{\rm Pic}}
\def\rg{{\rm rg}}
\def\et{{\rm et}}
\def\sep{{\rm s}}
\def\ppcm{{\rm ppcm}}
\def\ord{{\rm ord}}
\def\Gr{{\rm Gr}}
\def\ker{{\rm ker}}
\def\rk{{\rm rk}}
\def\Stab{{\rm Stab}}
\def\im{{\rm im}}
\def\Sm{{\rm Sm}}
\def\red{{\rm red}}
\def\Frob{{\rm Frob}}
\def\Ver{{\rm Ver}}


\def\beginProof{\par{\bf Proof. }}
 \def\endProof{${\qed}$\par\smallskip}
 \def\pr{^{\prime}}
 \def\prpr{^{\prime\prime}}
 \def\mtr#1{\overline{#1}}
 \def\ra{\rightarrow}
 \def\mfp{{\mathfrak p}}
 \def\mfm{{\mathfrak m}}
 
 \def\mQ{{\Bbb Q}}
 \def\mR{{\Bbb R}}
 \def\mZ{{\Bbb Z}}
 \def\mC{{\Bbb C}}
 \def\mN{{\Bbb N}}
 \def\mF{{\Bbb F}}
 \def\mA{{\Bbb A}}
  \def\mG{{\Bbb G}}
 \def\CI{{\cal I}}
 \def\CA{{\cal A}}
 \def\CE{{\cal E}}
 \def\CJ{{\cal J}}
 \def\CH{{\cal H}}
 \def\CO{{\cal O}}
 \def\CA{{\cal A}}
 \def\CB{{\cal B}}
 \def\CC{{\cal C}}
 \def\CK{{\cal K}}
 \def\CL{{\cal L}}
 \def\CI{{\cal I}}
 \def\CM{{\cal M}}
\def\CP{{\cal P}}
 \def\CZ{{\cal Z}}
\def\CR{{\cal R}}
\def\CG{{\cal G}}
\def\CX{{\cal X}}
\def\CY{{\cal Y}}
\def\CV{{\cal V}}
\def\CW{{\cal W}}
 \def\wt#1{\widetilde{#1}}
 \def\mod{{\rm mod\ }}
 \def\refeq#1{(\ref{#1})}
 \def\blb{{\big(}}
 \def\brb{{\big)}}
\def\mc{{{\mathfrak c}}}
\def\mcpr{{{\mathfrak c}'}}
\def\mcprpr{{{\mathfrak c}''}}
\def\ss{{\rm ss}}
\def\parf{{\rm parf}}
\def\P1{{{\bf P}^1}}
\def\cod{{\rm cod}}
\def\pr{\prime}
\def\prpr{\prime\prime}
\def\ss{\scriptstyle}
\def\OX{{ {\cal O}_X}}
\def\mpartial{{\mtr{\partial}}}
\def\inv{{\rm inv}}
\def\indlim{\underrightarrow{\lim}}
\def\prolim{\underleftarrow{\lim}}
\def\pprolim{'\prolim'}
\def\Pro{{\rm Pro}}
\def\Ind{{\rm Ind}}
\def\Ens{{\rm Ens}}
\def\without{\backslash}
\def\pbdb{{\Pro_b\ D^-_c}}
\def\qc{{\rm qc}}
\def\Com{{\rm Com}}
\def\an{{\rm an}}
\def\gfield{{\rm\bf k}}
\def\s{{\rm s}}
\def\dR{{\rm dR}}
\def\ari#1{\widehat{#1}}
\def\ul#1{\underline{#1}}
\def\sul#1{\underline{\scriptsize #1}}
\def\mou{{\mathfrak u}}
\def\ich{\mathfrak{ch}}
\def\cl{{\rm cl}}
\def\K{{\rm K}}
\def\R{{\rm R}}
\def\F{{\rm F}}
\def\L{{\rm L}}
\def\pgcd{{\rm pgcd}}
\def\rc{{\rm c}}
\def\N{{\rm N}}
\def\E{{\rm E}}
\def\H{{\rm H}}
\def\CHOW{{\rm CH}}
\def\A{{\rm A}}
\def\d{{\rm d}}
\def\Res{{\rm  Res}}
\def\GL{{\rm GL}}
\def\Alb{{\rm Alb}}
\def\alb{{\rm alb}}
\def\Hdg{{\rm Hdg}}
\def\Num{{\rm Num}}
\def\Irr{{\rm Irr}}
\def\Frac{{\rm Frac}}
\def\Sym{{\rm Sym}}
\def\TV{\rm TV}
\def\indlim{\underrightarrow{\lim}}
\def\prolim{\underleftarrow{\lim}}
\def\Hilb{{\rm Hilb}}
\def\Mor{{\rm Mor}}


\def\RHom{{\rm RHom}}
\def\rRHom{{\mathcal RHom}}
\def\rHom{{\mathcal Hom}}
\def\dotimes{{\overline{\otimes}}} 
\def\Ext{{\rm Ext}}
\def\rExt{{\mathcal Ext}}
\def\Tor{{\rm Tor}}
\def\rTor{{\mathcal Tor}}
\def\SP{{\mathfrak S}}

\def\H{{\rm H}}
\def\D{{\rm D}}
\def\Del{{\mathfrak D}}

\def\sh{{\rm sh}}

 \newtheorem{theor}{Theorem}[section]
 \newtheorem{prop}[theor]{Proposition}
 \newtheorem{propdef}[theor]{Proposition-Definition}
 \newtheorem{sublemma}[theor]{sublemma}
 \newtheorem{cor}[theor]{Corollary}
 \newtheorem{lemma}[theor]{Lemma}
 \newtheorem{sublem}[theor]{sub-lemma}
 \newtheorem{defin}[theor]{Definition}
 \newtheorem{conj}[theor]{Conjecture}
 \newtheorem{quest}[theor]{Question}

 \parindent=0pt
 \parskip=5pt

 \author{
 Damian R\"OSSLER\footnote{Institut de Math\'ematiques, 
Equipe Emile Picard, 
Universit\'e Paul Sabatier, 
118 Route de Narbonne,  
31062 Toulouse cedex 9, 
FRANCE, E-mail: rossler@math.univ-toulouse.fr}}
 \title{Infinitely $p$-divisible points on abelian varieties defined over function fields of characteristic $p>0$.}
\maketitle

\begin{abstract}
In this article we consider some questions raised by F. Benoist, E. Bouscaren and A. Pillay. 
We prove that infinitely $p$-divisible points on abelian varieties defined over function fields of transcendence degree one over 
a finite field are necessarily torsion points. We also prove that when the endomorphism ring of the abelian variety is $\mZ$ then there are no infinitely $p$-divisible points of order a power of $p$. 
\end{abstract}

\section{Introduction}

Fix once and for all a prime number $p$. 

 Let $K_0$ be the function 
field of a smooth curve over a finite field $\mF$ of characteristic $p$. Let $B$ be an abelian variety over $K_0$. 

For any abelian group $G$, define 
$
G^\#:=\cap_{l\geqslant 0}p^l G.
$ 
We call the elements of $G^\#$ the infinitely $p$-divisible points of $G$. 
Furthermore, write $\Tor(G)$ for the subset of $G$ consisting of elements of finite order. If $n\in\mN^*$, we write 
$\Tor^n(G)$ for the subset of $\Tor(G)$ consisting of elements of order prime to $n$. Similarly, we write $\Tor_n(G)$ for the subset of $\Tor(G)$ consisting of elements whose order divides a power of $n$. 

The following conjecture is made at the very end of the article  \cite{Benoist-Bouscaren-Pillay-Semiabelian}. 
\begin{conj}
The inclusion 
$
B(K^\sep_0)^\#\subseteq\Tor(B(K^\sep_0))
$ holds. 
\label{pilconj}
\end{conj}

Here $K_0^\sep$ is the separable closure of $K_0$. 

In the same context, F. Benoist then asked the following question (see \cite{Benoist-Communication}). 

\begin{quest}
Suppose that there are no non vanishing 
$\bar K_0$-homomorphisms $B_{\bar K_0}\to C$, where $C$ is an abelian variety, which has 
a model over $\bar\mF$. 

Do we have $B(K^\sep_0)^\#\cap \Tor_p(B(K_0^\sep))=0$ ?
\label{questb}
\end{quest}
We recall (in the notation of Question \ref{questb}) that a model of $C$ over $\bar\mF$ is an abelian variety 
$C_0$ over $\bar\mF$ such that $C\simeq C_0\times_{\bar\mF}\bar K_0$. 

In response to these questions, we shall prove the following two results 
in this text : 

\begin{theor}
Conjecture \ref{pilconj} holds. 
\label{maintheor}
\end{theor}

A. Pillay explained that Conjecture \ref{pilconj} can be viewed as a positive 
characteristic analog of Manin's "theorem of the kernel" (see 
\cite{Manin-Letter} for the latter). He also outlined a proof of the  Mordell-Lang conjecture over function fields 
of characteristic $p>0$ (in the case of abelian varieties), which is based on Conjecture \ref{pilconj} and a quantifier elimination result. 

Conjecture \ref{pilconj} is also linked to the Mordell-Lang conjecture in the following way. Lemma \ref{lemr} below (which plays a key role in the proof of Conjecture \ref{pilconj}) implies that the infinitely 
$p$-divisible points defined over a certain separable (but transcendental) field extension of $K_0$ are annulled by a fixed Weil polynomial (in particular it has no roots of unity among its roots) applied to a lifting of the Frobenius automorphism. This fact, combined with the existence of arc schemes (see for instance \cite[before Lemma 2.3]{Rossler-On-the-Manin})) as well as Proposition 6.1 in \cite{PR2} can be used to give a quick proof of Theorem 4.1 in \cite{Rossler-On-the-Manin}. This last theorem is the main tool in the proof of the Mordell-Lang conjecture given in 
\cite{Rossler-On-the-Manin}; more precisely, the Mordell-Lang conjecture follows quickly from it (see par. 3.2 in 
\cite{Rossler-On-the-Manin} for the argument), once the existence of jet schemes (see for instance 
par. 2.1 in \cite{Rossler-On-the-Manin} for the latter) is established.

Details about the argument outlined in the last paragraph will appear in \cite{Corpet1}, where it will also be generalized to the semiabelian situation; it would be very interesting to understand the link 
(if there is one) between this argument and A. Pillay's approach to the Mordell-Lang conjecture mentioned above. 

\begin{theor}
Suppose that $\End_{\bar K_0}(B)=\mZ$. Then  
{\it $B(K^\sep_0)^\#\cap \Tor_p(B(K_0^\sep))=0$.}
\label{bprop}
\end{theor}
Ie the answer to F. Benoist's question is affirmative if 
$\End_{\bar K_0}(B)=\mZ$. 

We actually prove a stronger result, but the stronger form is (as very often) more difficult to formulate. See the last remark of the text for this stronger form.

In particular, the answer to Question \ref{questb} is affirmative if $B$ is an 
elliptic curve over $K_0$. This was also proved in \cite{Benoist-Bouscaren-Pillay-Semiabelian}. 
In \cite[after Th. 3]{Voloch-Diophantine}, F. Voloch shows that if $B/K_0$ is ordinary and has maximal 
Kodaira-Spencer rank, then $B[p](K_0^\sep)=0$ and thus 
the answer to Question \ref{questb} is also affirmative in that situation. 

The hypothesis $B(K^\sep_0)^\#\cap \Tor_p(B(K_0^\sep))=0$ is a crucial hypothesis in B. Poonen's and J. F. Voloch's work 
on the Brauer-Manin obstruction over function fields. See for instance \cite[Th. B.]{Poonen-Voloch-The-Brauer-Manin}. 

The ideas underlying the proof of Theorem \ref{maintheor} are the following. First we notice that 
we may choose a (necessarily discrete) valuation on $K_0$ such that all the $p^l$-th roots of 
a given infinitely $p$-divisible point lie in a corresponding maximal unramified extension of 
the completion of $K_0$ along the valuation (this is Proposition \ref{prop1}). Secondly, we notice that an infinitely 
$p$-divisible point on $B$ can be recovered from the images of all its $p^l$-th roots in the reduction of $B$ modulo the 
ideal of the valuation (this is Lemma \ref{lemr}). This correspondence is Galois equivariant and this implies 
that an infinitely $p$-divisible point must have an infinite Galois orbit, if it is not a torsion 
point. This is a contradiction, because we are dealing with algebraic points (see the end of the proof of Proposition \ref{prop2}).

Here is the train of thought underlying the proof of Theorem \ref{bprop}. 
To obtain a contradiction, we suppose that we are given a sequence of 
non-zero points $P_i\in B(K_0^\sep)$ such that 
$p\cdot P_i=P_{i-1}$.   
We then consider the successive quotients of $B$ by the groups generated 
by the Galois orbits of the $P_i$. This gives a sequence of abelian varieties, 
which is shown to have decreasing modular height (see \refeq{brr}). Using a fundamental 
result of Zarhin (see Theorem \ref{zarth}), we deduce that infinitely many of the quotients 
are isomorphic and thus an abelian variety isogenous to $B$ is 
endowed with a non-trivial endomorphism. This contradicts one of the assumptions 
of Theorem \ref{bprop}. 

In the following two sections, we prove the results described 
above in the same order. The two sections are technically and terminologically 
 independent of each other and can be read in any order. 
 
{\bf Acknowledgments.} We are very grateful to R. Pink for his detailed comments 
on an earlier version of this text. Without his input, this text would certainly be a lot less clear. 
We also thank A. Pillay and F. Benoist for explaining their ideas and conjectures to us 
and for interesting discussions. Finally, we want to thank J.-F. Voloch for an interesting exchange on and around the problems addressed in this article. Thanks also to B. Poonen for his remarks. 

\section{Proof of Theorem \ref{maintheor}}

Before describing the proof, we would like to point out the following special case (which is not needed in the proof) : 

\begin{lemma}
If the $p$-rank of $B$ is $0$ then we have $
B(K^\sep_0)^\#\subseteq\Tor(B(K^\sep_0)).$
\end{lemma}
{\bf Proof} (of the lemma) This is a consequence of the Lang-N\'eron theorem. The details are 
left to the reader. 
 \endProof

Theorem \ref{maintheor} will be shown to be a simple consequence of the following two propositions. 

To formulate them, let $S$ be the spectrum of a complete noetherian local ring of characteristic $p$. Let $K$ be the fraction field of $S$. Let $\CA$ be an abelian scheme over $S$.   

Write $S^\sh$ for the spectrum of the strict henselisation of the ring underlying $S$ and $L$ for the fraction field of $S^\sh$. 
See for instance \cite[I, §4, p. 38]{Milne-Etale} for the definition of the strict henselization. 
We choose once and for all a $K$-embedding $L\hookrightarrow K^\sep$. 

\begin{prop}
Suppose that $S$ is a scheme of characteristic $p$. 
Suppose that the $p$-rank of the fibres of $\CA$ is constant. Then we have $$(p\cdot \CA(K^\sep))\cap \CA(L)=p\cdot \CA(L).$$
\label{prop1}
\end{prop}

The proof of Proposition \ref{prop1} will be given further below. 

\begin{prop}
Suppose that the residue field of the closed point of $S$ is finite. Then 
we have $A(L)^\#=\Tor(A(L)^\#)$. 
\label{prop2}
\end{prop}

The proof of Proposition \ref{prop2} will be given further below.

{\bf Proof of Theorem \ref{maintheor}} (assuming Propositions \ref{prop1} \& \ref{prop2}). Let $x\in B(K_0^\sep)^\#$. We shall show that 
$x$ is a torsion point. 

We may assume without restriction 
of generality that $x\in B(K_0)$ (otherwise, we replace $K_0$ by $K_0(\kappa(x))$). After again possibly replacing $K_0$ by a finite extension, we may suppose that $B[p](K_0)=B[p](\bar K_0)$. 
Let $U$ be a smooth algebraic curve over $\mF$. 
We may suppose without restriction of generality that there is an abelian scheme $\CB$ over $U$ extending $B$ and that 
the $p$-rank of the fibres of $\CB$ is constant. Let now $u$ be any closed point of $U$. Let $S$ be the completion of the local 
ring of $U$ at $u$ and let $\CA:=\CB_{S}$. Let $K$ be the fraction field of $S$ and $L$ be the fraction field of 
the strict henselisation $S^\sh$ of $S$, as before. 

Now let $x_1, x_2,\dots\in B(K_0^\sep)$ be elements such that $p\cdot x_l=x_{l-1}$ for all $l\geqslant 2$ and such that 
$p\cdot x_1=x$. Such a sequence exists by assumption. Applying Proposition \ref{prop1} to $x$ and all the $x_l$ successively, 
we conclude that $x\in \CA(L)^\#$. We conclude the proof of Theorem \ref{maintheor} by appealing to Proposition \ref{prop2}.\endProof

{\bf Proof of Proposition \ref{prop1}}.
Let $N:=\ker\ [p]$ be the kernel of the multiplication by $p$ morphism $[p]:\CA\to\CA$. 
This is a (non-reduced) finite group-scheme over $S$. There is an exact sequence of finite, flat, 
$S$-group schemes
$$
0\to N_0\to N\to N_\et\to 0\ \ \ \ \ \ (\%)
$$
where $N_0$ is connected and $N_\et$ is \'etale. See \cite[(1.4)]{Tate-p-divisible} for this. Since the unit 
section of $N$ is open and closed (by the assumption on the $p$-rank), we see that $N_0$ is an infinitesimal group-scheme.

Thus the decomposition (\%) of $N$ leads to $S$-isogenies of abelian schemes
$$
\CA\stackrel{\phi}{\to} \CA_{1}\stackrel{\mu}{\to} \CA
$$
such that $\mu$ is \'etale, $\phi$ is purely inseparable and $\mu\circ\phi=[p]$. 

Let now $x\in (p\cdot A(K^\sep))\cap A(L)$. 
Let $T\hookrightarrow\CA_{S^\sh}$ be the closed 
subscheme arising from the section of  $\CA_{S^\sh}$ over $S^\sh$ associated to $x$. We have $S^\sh$-morphisms 
$[p]^*T\to \mu^*_{S^\sh}T\to T$, where the morphism $[p]^*T\to \mu^*_{S^\sh}T$ is purely inseparable over $L$ and 
the morphism $\mu^*_{S^\sh}T\to T$ is \'etale. Furthermore, by construction, there is an element $y\in \CA(K^{\sep})$ such 
that $y\in ([p]^*T)(K^{\sep})$ and such that $p\cdot y=x$. Let $V$ be the connected component of 
$\mu^*T$ containing the image of $\phi_{S^\sh}(y)$. We obtain a sequence of $S^\sh$-morphisms 
$$
\Spec\ K^{\sep}\to ([p]^*T)_V\to V\to T\to A
$$
whose composition is $x$.

Since $S^\sh$ is strictly henselian, the morphism $V\to T$ must be an isomorphism.
Let now $v\in V$ be the image of the composition $\Spec\ K^{\sep}\to ([p]^*T)_V\to V$. Let 
$w\in ([p]^*T)_V$ be the image of $y$. Since the morphism $([p]^*T)_V\to V$ is purely inseparable, 
the extension of residue fields $\kappa(w)|\kappa(v)$ is purely inseparable. Hence it must be trivial, since 
the field extension $\kappa(w)|\kappa(T)=L=\kappa(v)$ is separable, because it is a subextension of
 the separable algebraic extension $K^{\sep}|\kappa(T)$. Hence 
the extension $\kappa(w)|\kappa(T)$ is trivial. 

This shows that we actually have $x\in p\cdot A(L)$, which completes the proof of the proposition.
\endProof

{\bf Proof of Proposition \ref{prop2}}. Let us temporarily use the following notation. Let $R:=\Gamma(S,\CO_S)$ be the local ring underlying $S$. 
Let $\mfm\subseteq R$ be the maximal ideal of $R$ and let $k:=R/\mfm$. We write correspondingly 
$R^\sh:=\Gamma(S^\sh,\CO_{S^\sh})$ and $\mfm^\sh\subseteq R^\sh$ for the maximal ideal of $R^\sh$. 
We fix an identification of $\bar k$ with $R^\sh/\mfm^\sh$. Let also $\widehat{S^\sh}$ be the spectrum of the completion $\widehat{R^\sh}$ of $R^\sh$ along $\mfm^\sh$ and let 
$\widehat{L}$ be the fraction field of $\widehat{R^\sh}$. 

Let $A_0$ be the special fibre of $\CA$, ie the fibre of $\CA$ over the unique closed point of $S$. 

Let ${I}_p(A_0(\bar k )):=\varprojlim_{l\geqslant 0}A_0(\bar k )$, where the transition morphisms are all 
given by multiplication by $p$ (beware : $I_p(A_0(\bar k ))$ is not the classical Tate module). So 
an element of $\varprojlim_{l\geqslant 0}A_0(\bar k )$ is given by an inverse system
$$
\cdots\to t_{l}\to t_{l-1}\to\cdots\to t_0
$$
where $t_l\in A_0(\bar k)$ and $t_{l-1}=p\cdot t_l$ for all $l\geqslant 0$. Now notice that for any $l\geqslant 1$, the kernel of the reduction map $\CA(R^\sh/\mfm^{\sh,l})\to A_0(\bar k)$ is a commutative group, which is killed by multiplication by $p^{l-1}$. 
This follows from the fact that the kernel of the reduction map $\CA(\widehat{S^\sh})\to A_0(\bar k)$ can be identified 
with the points of a commutative formal group. 
See \cite[Th. 8.5.9 (a), p. 213]{Fantechi-FGA}.

For a fixed $r\geqslant 1$, let $\lambda_{r}:I_p(A_0(\bar k ))\to \CA(R^\sh/\mfm^{\sh,r})$ be the map sending the 
inverse system
$$
\cdots\to t_{r}\to t_{r-1}\to\cdots\to t_0
$$ 
to the element $p^{r-1}\cdot \widetilde{t}_{r-1}$, where 
$\widetilde{t}_{r-1}$ is any lifting of $t_{r-1}$ to  $\CA(R^\sh/\mfm^{\sh,r})$. 
Notice that the composition of $\lambda_{r}$ with the reduction map 
$\CA(R^\sh/\mfm^{\sh,r})\to \CA(R^\sh/\mfm^{\sh,r-1})$ is the map $\lambda_{r-1}$. Hence we obtain 
a homomorphism 
$$
\lambda:I_p(A_0(\bar k ))\to \varprojlim_r \CA(R^\sh/\mfm^{\sh,r})\simeq\CA(\widehat{S^\sh})\simeq \CA(\widehat{L}).
$$ 

{\bf Remark.} A variant of the map $\lambda$ appears in the theory of the Serre-Tate lifting of ordinary abelian varieties 
over finite fields (see \cite{Katz-Serre-Tate}). 

\begin{lemma}
We have $\CA(\widehat{L})^\#={\rm Im}(\lambda)$.
\label{lemr}
\end{lemma}

Lemma \ref{lemr} is a variant of \cite[Lemme 3.2.1]{Raynaud-Sous-var}. 
The proofs of both Lemmata are actually identical. 

{\bf Proof} (of the lemma). 
Let $x_0, x_1,\dots \in A(\widehat{L})$ be such that 
$p\cdot x_i=x_{i-1}$ for all $i\geqslant 1$. Let $\rho: A(\widehat{L})\to A_0(\bar k )$ be the reduction map. 
Consider the inverse system $q\in I_p(A_0(\bar k ))$ given by $$\cdots\to \rho(x_{i})\to \rho(x_{i-1})\to\cdots\to\rho(x_0).$$ We claim that $\lambda(q)=x_0$. For this, it is sufficient to prove that 
$$
\lambda_r(q)=x_0\ (\mod\ \mfm^{\sh,r})
$$
for all $r\geqslant 1$. We compute
$$
\lambda_r(q)=p^{r-1}\cdot (x_{r-1} (\mod\ \mfm^{\sh,r}))=(p^{r-1}\cdot x_{r-1}) (\mod\ \mfm^{\sh,r}))= 
x_0\ (\mod\ \mfm^{\sh,r}).
$$
This shows that $\cap_{l\geqslant 0}p^l A(\widehat{L})\subseteq{\rm Im}(\lambda)$. 

For the opposite inclusion (which is actually not needed in the text), let 
$$
\cdots\to t_1\to t_0
$$
be an element of $I_p(A_0(\bar k ))$. This element is the 
image under multiplication by $p$ of the element
$$
\cdots\to t_2\to t_1
$$
(\ in other words $I_p(A_0(\bar k ))^\#=I_p(A_0(\bar k ))$\ ). Hence 
any image by $\lambda$ of an element of $I_p(A_0(\bar k ))$ is infinitely $p$-divisible.
\endProof

{\bf Remark.} Until now in the proof, we only used the fact that $k$ is a field of characteristic $p$ 
(not that it is a finite field). In particular, the last lemma is true in that generality. 

Now notice that the map $\lambda$ is compatible with the natural action of $\Gal(\bar k |k)$ 
on $I_p(A_0(\bar k ))$ and $\CA(\widehat{L})$. Let 
$\sigma\in\Gal(\bar k |k)\simeq\widehat{\mZ}$ be the Frobenius automorphism. 
Let $x\in \CA(L)^\#\subseteq\CA(\widehat{L})^\#={\rm Im}(\lambda)$. 
Since the residue field of $x$ lies in an algebraic extension of $K$, there exists $r\geqslant 1$ such that 
$\sigma^r(x)=x$. On the other hand, by the Weil conjectures, there is a polynomial $Q$ with integer coefficients and no roots of unity among its roots, such that  $Q(\sigma^r)=0$ on $I_p(A_0(\bar k ))$. Hence $Q(\sigma^r(x))=Q(1)\cdot x=0$. 
Since $Q(1)\not= 0$, this implies that $x$ is a torsion point. This proves Proposition \ref{prop2}. 
\endProof

{\bf Remark.} Let $L_1$ be the algebraic closure of $L$ in $\widehat{L}$ 
(notice that if $S$ is a discrete valuation ring then $L_1=L$; see for instance \cite[Ex. 8.3.34, p. 360]{Liu-Algebraic}). The above proof of Proposition \ref{prop2} actually shows that 
$$
\CA(\widehat{L})^\#\cap\CA(L_1)=\Tor(\CA(L_1)),
$$
which is stronger. 

\section{Proof of Theorem \ref{bprop}}

This section can be read independently of the previous one. 

To prove Theorem \ref{bprop}, we shall crucially need the following theorem. 
Let $U$ be a smooth proper curve over $\mF$, whose function field is $K_0$. 
A semiabelian scheme $\CA/U$ is a commutative group scheme over $U$, whose 
geometric fibres are semiabelian varieties. We shall say that a semiabelian scheme over 
$U$ is generically abelian if its geometric generic fibre is an abelian variety. 
If $\CA/U$ is a semiabelian scheme with zero-section $\epsilon$, we shall 
write
$$
\omega_{\CA/U}^0:=\epsilon^*(\det(\Omega_{\CA/U}))
$$
for the determinant of the sheaf of differentials of $\CA$ over $U$, restricted to 
the zero-section. 

%

\begin{theor}[Zarhin] 
Let $\beta,g,n\geqslant 0$. Suppose that $(n,p)=1$ and that 
$n\geqslant 3$. 
Up to isomorphism, there is only a finite number of generically abelian semiabelian schemes $\CA/U$ of relative 
dimension $g$ over $U$, such that 
\begin{itemize}
\item $\deg_{U}(\omega^0_{\CA/U})=\beta$;
\item there exists an isomorphism of group schemes 
$\CA_{K_0}[n]\simeq(\mZ/n\mZ)^{2g}_{K_0}$.
\end{itemize}
\label{zarth}
\end{theor}

The Theorem \ref{zarth} is well-known to specialists but I could find no formal proof of it in the literature. 

{\bf Proof}  (of Theorem \ref{zarth}). 
We shall first prove the following statement. Let $\beta\in\mN$. 

Up to isomorphism, there is only a finite number of generically abelian semiabelian schemes $\CA/U$ of relative 
dimension $g$ over $U$, such that (*)
\begin{itemize}
\item there exists a principal polarisation on $\CA_{K_0}$;
\item there exists a symplectic isomorphism $(\mZ/n\mZ)^{2g}_{K_0}\simeq\CA_{K_0}[n]$; 
\item $\deg_{U}(\omega^0_{\CA/U})=\beta$. 
\end{itemize}

To prove this, we shall use the following deep results of D. Mumford, A. Grothendieck, L. Moret-Bailly, C.-L. Chai and G. Faltings. 

Let ${\bf A}_{g,n}$ be the functor  from the category of locally noetherian $\mF_p$-schemes to the category of 
sets, such that
\begin{eqnarray*}
&&{\bf A}_{g,n}(S)=\{\,\textrm{\rm isomorphism classes of the following objects :}\\
&&\textrm{principally polarized abelian schemes over $S$ endowed}\\
&&\textrm{with 
a symplectic isomorphism $(\mZ/n\mZ)^{2g}_{S}\simeq\CA[n]$}\,\}
\end{eqnarray*}

Then D. Mumford proves (see \cite{Mumford-GIT}) that  the functor ${\bf A}_{g,n}$ is representable 
by a scheme, which is separated and of finite type over $\mF_p$. We shall also denote 
this scheme by  ${\bf A}_{g,n}$. 

Furthermore, in \cite[V, 2., Th. 2.5]{Faltings-Chai-Degeneration}, C. Chai and G. Faltings prove that there exists
\begin{itemize}
\item a scheme ${\bf A}_{g,n}^*$, which is proper over $\mF_p$;
\item an open immersion ${\bf A}_{g,n}\hookrightarrow {\bf A}_{g,n}^*$;
\item a semiabelian scheme $\CB$ over ${\bf A}_{g,n}^*$, such that $\CB_{{\bf A}_{g,n}}$ is isomorphic 
to the abelian scheme underlying the universal object over ${\bf A}_{g,n}$.
\end{itemize}
Also they show that $\omega^0_{\CB/{\bf A}_{g,n}^*}$ is an ample line bundle. 

Now write $Z:=U\times_{\mF_p}{\bf A}_{g,n}^*$. Recall that the Hilbert scheme $\Hilb(Z/\mF_p)$ 
is a scheme over $\mF_p$, which is locally of finite type and such that 
$$
\Hilb(Z/\mF_p)(T)=\{\textrm{closed subschemes of $Z_T$, which are 
proper and flat over $T$}\}
$$
for any locally noetherian scheme $T$ over $\mF_p$ (see \cite{FGA-221}). 

Furthermore, fix $\Phi\in\mQ[\lambda]$, a polynomial 
with rational coefficients and $L/Z$ an ample line bundle. The $\mF_p$-scheme $\Hilb_\Phi(Z/\mF_p)$ is then 
characterized by the property that 
\begin{eqnarray*}
&&\Hilb_\Phi(Z/\mF_p)(T):=\\
&&\{\textrm{closed subschemes $W$ of $Z_T$, which are 
proper and flat over $T$}\\
&&\textrm{and such that $\chi(W_t,L_{W_t}^{\otimes\lambda})=\Phi(\lambda)$ for all $\lambda\in\mN$ and all $t\in T$}\}\\
\end{eqnarray*}
Here $W_t$ is the fibre at $t\in T$ of the morphism $W\to T$ and $L_{W_t}$ is the pull-back of $L$ to 
$W_t$ by the natural morphism $W_t\to Z$. The symbol $\chi(\cdot)$ refers to the Euler characteristic; 
by definition
$$
\chi(W_t,L_{W_t}^{\otimes\lambda})=\sum_{r\geqslant 0}(-1)^r\dim_{\kappa(t)}H^r(W_t,L_{W_t}^{\otimes\lambda}).
$$
It is shown in \cite{FGA-221}, that $\Hilb_\Phi(Z/\mF_p)$ is projective over $\mF_p$ (as a consequence of the 
projectivity of $Z$). Notice that by construction, we have a decomposition
$$
\Hilb(Z/\mF_p)=\coprod_{\Phi\in\mQ[\lambda]}\Hilb_\Phi(Z/\mF_p).
$$
Finally, it is shown in 
\cite[part II, 5.23]{Fantechi-FGA} that the functor $\Mor_{\mF_p}(U,{\bf A}^*_{g,n})$ from 
locally noetherian $\mF_p$-schemes $T$ to the category of sets, such that 
$$
\Mor_{\mF_p}(U,{\bf A}^*_{g,n})(T)=\{\textrm{$T$-morphisms from $U_T$ to ${\bf A}^*_{g,n,T}$}\}
$$
is representable by an open subscheme of $\Hilb(Z/\mF_p)$. The open immersion 
$\Mor_{\mF_p}(U,{\bf A}^*_{g,n})\hookrightarrow\Hilb(U\times_{\mF_p}{\bf A}_{g,n}^*/\mF_p)$ is described 
by the natural transformation of functors
$$
\textrm{$T$-morphism $f$ from $U_T$ to ${\bf A}^*_{g,n,T}$}\mapsto\textrm{graph of $f$}
$$

Now fix an ample line bundle $M$ on $U$. Let $L$ be the line bundle $M\boxtimes\omega^0_{\CB/{\bf A}^*_{g,n}}$ 
on $Z=U\times_{\mF_p}{\bf A}_{g,n}^*$. 

Now we are finally ready to tackle our proof of finiteness. Suppose that we are given $\CA/U$ as in (*). Restricting to the generic point of $U$, we get 
a morphism $\Spec\,\kappa(U)\to {\bf A}_{g,n}^*$ (whose image is in ${\bf A}_{g,n}$) and since 
${\bf A}_{g,n}^*$ is proper over $\mF_p$ and $U$ is a Dedekind scheme, this extends to a morphism 
$\phi:U\to  {\bf A}_{g,n}^*$. We thus get a point $\phi\in\Mor_{\mF_p}(U,{\bf A}^*_{g,n})(\mF_p)$. 
Let $\Gamma_\phi\hookrightarrow U\times_{\mF_p}{\bf A}^*_{g,n}$ be the graph of $\phi$. We 
compute its Hilbert polynomial
\begin{eqnarray}
\chi(\Gamma_\phi,L_\Gamma^{\otimes\lambda})&=&\chi(U,(M\otimes\phi^*(\omega^0_{\CB/{\bf A}^*_{g,n}}))^{\otimes\lambda})=\deg_U((M\otimes\phi^*(\omega^0_{\CB/{\bf A}^*_{g,n}}))^{\otimes\lambda})+1-g_U\nonumber\\&=&
\lambda\cdot\deg_U(M\otimes\phi^*(\omega^0_{\CB/{\bf A}^*_{g,n}}))+1-g_U\nonumber\\
&=&\lambda\cdot\deg_U(M\otimes\omega_{\CA/U})+1-g_U\nonumber\\
&=&
\lambda\cdot\deg_U(M)+\lambda\cdot\deg_U(\omega_{\CA/U})+1-g_U\nonumber\\
&=&
\lambda\cdot\deg_U(M)+\lambda\cdot\beta+1-g_U:=Q(\lambda).
\label{creq}
\end{eqnarray}
Here $g_U$ is the genus of $U$. The second equality is justified by the Riemann-Roch theorem on $U$. 
From \refeq{creq} we see that the Hilbert polynomial of $\phi\in\Mor_{\mF_p}(U,{\bf A}^*_{g,n})(\mF_p)$ 
only depends on $\beta$ (once $M$ is given) and thus\linebreak  
\mbox{$\phi\in \Hilb_{Q(\lambda)}(U\times_{\mF_p}{\bf A}_{g,n}^*/\mF_p)(\mF_p)$.} 

Now to prove that there are only a finite number of generically abelian semiabelian scheme 
$\CA/U$ satisfying (*), just notice that the set $\Hilb_{Q(\lambda)}(U\times_{\mF_p}{\bf A}_{g,n}^*/\mF_p)(\mF_p)$ 
is finite, since $\Hilb_{Q(\lambda)}(U\times_{\mF_p}{\bf A}_{g,n}^*/\mF_p)$ is projective and hence of finite type 
over $\mF_p$. 

To conclude, recall the following facts.

First, a generically abelian semiabelian scheme $\CA/U$ is completely determined by its generic fibre $\CA_{K_0}$ - see \cite[IX, cor. 1.4, p. 130]{Raynaud-Faisceaux}. Secondly, by Zarhin's "trick" (see \cite[IX, 1.]{Moret-Bailly-Pinceaux}), 
for any abelian variety $C/K_0$, the abelian variety 
$(C\times_{K_0}C^\vee)^4$ can be principally polarized. Thirdly, 
a given abelian variety $C/K_0$ has only a finite number of 
direct factors (see \cite[Th. 18.7]{Milne-Abelian}). 
Fourthly, if $C/K_0$ is an abelian variety, 
which extends to a semiabelian scheme $\CC$ over $U$, then the dual 
abelian variety $C^\vee$ has the same property (by Grothendieck's criterion 
\cite[IX, Prop. 5.13 (c)]{SGA7.1}) and furthermore 
$\deg(\omega_{\CC/U}^0)=\deg(\omega^0_{\CC^\vee/U})$ 
(here we wrote somewhat sloppily $\CC^\vee$ for the semiabelian extension 
of $C^\vee$). See \cite[V, 3., Lemma 3.4]{Faltings-Chai-Degeneration} for the latter equality. 

Putting all these facts together, the theorem readily follows from the just proven fact that 
there are only a finite number of generically abelian semiabelian scheme 
$\CA/U$ satisfying (*).\endProof

{\bf Proof}  (of Theorem \ref{bprop}). By contradiction. So suppose 
that there exists a point $P\in B(K^\sep_0)^\#$ such that $P\not=0$ and 
$p^r\cdot P=0$ for some $r\geqslant 1$. Then there exists a sequence of 
points $(P_n\in B(K^\sep_0))_{n\in\mN^*}$ such that 
$P=P_1$ and $p\cdot P_{i+1}=P_i$ for all $i\in\mN^*$. Let 
$G_i$ be the subgroup of $B(K^\sep_0)^\#$ generated by the elements 
$\gamma(P_i)$, where $\gamma$ runs through $\Gal(K_0^\sep|K_0)$. 
By Galois descent, the groups $G_i$ are naturally defined over $K_0$ 
and we have natural inclusion morphisms $G_1\hookrightarrow 
G_2\hookrightarrow\dots$ over $K_0$. Furthermore, the order $d_i$ of the group 
$G_i$ is strictly increasing as a function of $i\in\mN^*$. 

Now we may replace $K_0$ by a finite extension and $U$ by the corresponding 
projective curve over $k$, without restricting generality. 
Hence we may suppose that $B$ extends to a generically abelian 
semiabelian scheme $\CB$ over $U$ (by Grothendieck's semiabelian reduction theorem 
\cite[IX]{SGA7.1}). 
Furthermore, we may suppose that there is an isomorphism 
$\rho_B:(\mZ/l\mZ)^{2g}_{K_0}\simeq\CB_{K_0}[l]$, where $l\geqslant 3$ is a 
natural number prime to $p$. We fix such an $l$ for the rest of the proof. 

Now look at the sequence of $K_0$-isogenies
\begin{eqnarray}
B\to B/G_1\to B/(G_1\cdot G_2)\to\cdots B/(G_1\cdot G_2\cdot G_3)\to\dots
\label{brr}
\end{eqnarray}
Let $B_i:=B/(G_1\cdots G_i)$ and $B_0:=B$. Let $\pi_i:\CB_i\to U$ be the connected component 
of the N\'eron model of $B_i$ over $U$. Since 
all the $B_i$ are isogenous, the criterion \cite[IX, Prop. 5.13 (c)]{SGA7.1} shows that 
all the $\CB_i$ are generically abelian semiabelian schemes over $U$. Furthermore, 
by the universal property of N\'eron models, the morphisms in \refeq{brr} extend 
to $U$-morphisms
\begin{eqnarray}
\CB\stackrel{\phi_0}{\to}\CB_1\stackrel{\phi_1}{\to}\CB_2\to\cdots.
\label{brr1}
\end{eqnarray}
Finally, since the isogenies $\phi_i$ are of degree prime to $l$, 
we obtain isomorphisms $\phi_{i}\circ\phi_{i-1}\circ\cdots\phi_0\circ\rho_B:(\mZ/l\mZ)^{2g}_{K_0}\simeq\CB_{i,K_0}[l]$  for all $i$. 

By construction, the morphisms $\phi_i$ are generically \'etale. 
Looking at differentials, we see that there is an exact sequence of coherent sheaves
$$
\phi_i^*\Omega_{\CB_{i+1}/U} \to\Omega_{\CB_{i}/U}\to\Omega_{\phi_i}\to 0
$$
since $\Omega_{\CB_{i+1}/U}$ is locally free and since the morphism of sheaves 
$\phi_i^*\Omega_{\CB_{i+1}/U} \to\Omega_{\CB_{i}/U}$ is injective at the generic point of 
$\CB_i$, we see that the morphism \mbox{$\phi_i^*\Omega_{\CB_{i+1}/U} \to\Omega_{\CB_{i}/U}$} is actually 
injective so that we have an exact sequence
\begin{equation}
0\to\phi_i^*\Omega_{\CB_{i+1}/U} \to\Omega_{\CB_{i}/U}\to\Omega_{\phi_i}\to 0.
\label{crapeq}
\end{equation}
Restricting the sequence \refeq{crapeq} to the zero section $\epsilon_i:U\to\CB_i$, we get
\begin{equation}
0\to\epsilon^*_i(\Omega_{\CB_{i+1}/U})\to \epsilon^*_i(\Omega_{\CB_{i}/U})\to\epsilon^*_i(\Omega_{\phi_i})\to 0.
\label{seqint}
\end{equation}
Here the exactness on the left is again justified by the fact that the sequence \refeq{seqint} 
is generically exact on $U$ and the fact that $\epsilon^*_{i+1}(\Omega_{\CB_{i+1}/U})$ 
is locally free. 
Applying $\deg_U(\cdot)$ to the objects in the last sequence, we conclude that 
$$
\deg_{U}(\omega^0_{\CB_i/U})-\deg_U(\epsilon^*_{i}(\Omega_{\phi_i}))=\deg_{U}(\omega^0_{\CB_{i+1}/U})
$$
Since $\epsilon^*_{i}(\Omega_{\phi_i})$ is a torsion sheaf, we have $\deg_U(\epsilon^*_i(\Omega_{\phi_i}))\geqslant 0$ 
(see \cite[Ex. 6.12, chap. 7, p. 149]{Hartshorne-Algebraic}).

So we conclude that $\deg_{U}(\omega^0_{\CB_i/U})\leqslant \deg_{U}(\omega^0_{\CB/U})$ for all 
$i\in\mN$. From Theorem \ref{zarth}, we conclude that there exists $j,l\in\mN$ such that $j<l$ 
and such that $\CB_j\simeq\CB_l$. Since there is also an \'etale $K_0$-isogeny 
$B_j\to B_l$ of degree divisible by $p^{l-j}$, this means that there is an \'etale $K_0$-isogeny 
$B_j\to B_j$ of degree divisible by $p^{l-j}$. Since this last isogeny cannot coincide with 
the isogeny given by a power of $p$ (because the isogeny given by multiplication by $p$ always has 
an inseparable part), this implies that $\End_{\bar K_0}(B_j)_\mQ\not=\mQ$. Thus $\End_{\bar K_0}(B)_\mQ\not=\mQ$. 
This concludes the proof.\endProof

{\bf Remark.} The Proof of Theorem \ref{bprop} gives a somewhat stronger result. 
In fact it shows that the conclusion of Theorem \ref{bprop} holds if the following weaker assumption holds : 
there is no $\bar K_0$-isogeny $\phi:B_{\bar K_0}\to B'$, such that 
\begin{itemize}
\item $\phi$ is \'etale; 
\item $\deg(\phi)$ is a power of $p$;
\item $B'$ carries an \'etale and finite $\bar K_0$-endomorphism, whose degree is a power of $p$. 
\end{itemize}


\begin{bibdiv}
\begin{biblist}

\bib{Voloch-Towards}{article}{
  author={Abramovich, Dan},
  author={Voloch, Jos{\'e} Felipe},
  title={Toward a proof of the Mordell-Lang conjecture in characteristic $p$},
  journal={Internat. Math. Res. Notices},
  date={1992},
  number={5},
  pages={103--115},
  issn={1073-7928},
}

\bib{Benoist-Communication}{article}{
  author={Benoist, Franck},
  title={A theorem of the Kernel in characteristic p},
  status={Preprint. February 8th 2011.},
}

\bib{Benoist-Bouscaren-Pillay-Semiabelian}{article}{
  author={Benoist, Franck},
  author={Bouscaren, Elisabeth},
  author={ Pillay, Anand},
  title={Semiabelian varieties over separably closed fields, maximal divisible subgroups, and exact sequences},
  status={Preprint arXiv:0904.2083v1},
}

\bib{Corpet1}{article}{
  author={C. Corpet},
  title={Around the Mordell-Lang and Manin-Mumford conjectures in positive characteristic},
  status={in preparation},
}

\bib{Faltings-Chai-Degeneration}{book}{
  author={Faltings, Gerd},
  author={Chai, Ching-Li},
  title={Degeneration of abelian varieties},
  series={Ergebnisse der Mathematik und ihrer Grenzgebiete (3) [Results in Mathematics and Related Areas (3)]},
  volume={22},
  note={With an appendix by David Mumford},
  publisher={Springer-Verlag},
  place={Berlin},
  date={1990},
  pages={xii+316},
  isbn={3-540-52015-5},
}

\bib{Fantechi-FGA}{collection}{
  author={Fantechi, Barbara},
  author={G{\"o}ttsche, Lothar},
  author={Illusie, Luc},
  author={Kleiman, Steven L.},
  author={Nitsure, Nitin},
  author={Vistoli, Angelo},
  title={Fundamental algebraic geometry},
  series={Mathematical Surveys and Monographs},
  volume={123},
  note={Grothendieck's FGA explained},
  publisher={American Mathematical Society},
  place={Providence, RI},
  date={2005},
  pages={x+339},
  isbn={0-8218-3541-6},
}

\bib{FGA-221}{article}{
  author={Grothendieck, Alexander},
  title={Techniques de construction et th\'eor\`emes d'existence en g\'eom\'etrie alg\'ebrique. IV. Les sch\'emas de Hilbert},
  language={French},
  conference={ title={S\'eminaire Bourbaki, Vol.\ 6}, },
  book={ publisher={Soc. Math. France}, place={Paris}, },
  date={1995},
  pages={Exp.\ No.\ 221, 249--276},
}


\bib{SGA3-2}{book}{
  title={Sch\'emas en groupes. II: Groupes de type multiplicatif, et structure des sch\'emas en groupes g\'en\'eraux},
  language={},
  series={S\'eminaire de G\'eom\'etrie Alg\'ebrique du Bois Marie 1962/64 (SGA 3). Dirig\'e par M. Demazure et A. Grothendieck. Lecture Notes in Mathematics, Vol. 152},
  publisher={Springer-Verlag},
  place={Berlin},
  date={1962/1964},
  pages={ix+654},
}

\bib{SGA7.1}{book}{
  title={Groupes de monodromie en g\'eom\'etrie alg\'ebrique. I},
  language={},
  series={Lecture Notes in Mathematics, Vol. 288},
  note={S\'eminaire de G\'eom\'etrie Alg\'ebrique du Bois-Marie 1967--1969 (SGA 7 I); Dirig\'e par A. Grothendieck. Avec la collaboration de M. Raynaud et D. S. Rim},
  publisher={Springer-Verlag},
  place={Berlin},
  date={1972},
  pages={viii+523},
}

\bib{Hartshorne-Algebraic}{book}{
  author={Hartshorne, Robin},
  title={Algebraic geometry},
  note={Graduate Texts in Mathematics, No. 52},
  publisher={Springer-Verlag},
  place={New York},
  date={1977},
  pages={xvi+496},
  isbn={0-387-90244-9},
}

\bib{Hrushovski-Mordell-Lang}{article}{
  author={Hrushovski, Ehud},
  title={The Mordell-Lang conjecture for function fields},
  journal={J. Amer. Math. Soc.},
  volume={9},
  date={1996},
  number={3},
  pages={667--690},
  issn={0894-0347},
}

\bib{Katz-Serre-Tate}{article}{
  author={Katz, Nicholas},
  title={Serre-Tate local moduli},
  conference={ title={Algebraic surfaces}, address={Orsay}, date={1976--78}, },
  book={ series={Lecture Notes in Math.}, volume={868}, publisher={Springer}, place={Berlin}, },
  date={1981},
  pages={138--202},
}

\bib{Liu-Algebraic}{book}{
  author={Liu, Qing},
  title={Algebraic geometry and arithmetic curves},
  series={Oxford Graduate Texts in Mathematics},
  volume={6},
  note={Translated from the French by Reinie Ern\'e; Oxford Science Publications},
  publisher={Oxford University Press},
  place={Oxford},
  date={2002},
  pages={xvi+576},
  isbn={0-19-850284-2},
}

\bib{Manin-Letter}{article}{
  author={Manin, Yuri I.},
  title={Letter to the editors: ``Rational points on algebraic curves over function fields'' [Izv.\ Akad.\ Nauk SSSR Ser.\ Mat.\ {\bf 27} (1963), 1397--1442; MR0157971 (28 \#1199)]},
  language={Russian},
  journal={Izv. Akad. Nauk SSSR Ser. Mat.},
  volume={53},
  date={1989},
  number={2},
  pages={447--448},
  issn={0373-2436},
  translation={ journal={Math. USSR-Izv.}, volume={34}, date={1990}, number={2}, pages={465--466}, issn={0025-5726}, },
  review={\MR {998307 (90f:11039)}},
}

\bib{Milne-Abelian}{article}{
  author={Milne, James S.},
  title={Abelian varieties},
  conference={ title={Arithmetic geometry}, address={Storrs, Conn.}, date={1984}, },
  book={ publisher={Springer}, place={New York}, },
  date={1986},
  pages={103--150},
}

\bib{Milne-Etale}{book}{
   author={Milne, James S.},
   title={\'Etale cohomology},
   series={Princeton Mathematical Series},
   volume={33},
   publisher={Princeton University Press},
   place={Princeton, N.J.},
   date={1980},
   pages={xiii+323},
   isbn={0-691-08238-3},
}

\bib{Moret-Bailly-Pinceaux}{article}{
  author={Moret-Bailly, Laurent},
  title={Pinceaux de vari\'et\'es ab\'eliennes},
  language={},
  journal={Ast\'erisque},
  number={129},
  date={1985},
  pages={266},
  issn={0303-1179},
}

\bib{Mumford-GIT}{book}{
  author={Mumford, David},
  author={Fogarty, John},
  author={Kirwan, Frances},
  title={Geometric invariant theory},
  series={Ergebnisse der Mathematik und ihrer Grenzgebiete (2) [Results in Mathematics and Related Areas (2)]},
  volume={34},
  edition={3},
  publisher={Springer-Verlag},
  place={Berlin},
  date={1994},
  pages={xiv+292},
  isbn={3-540-56963-4},
}

\bib{PR2}{article}{
  author={Pink, Richard},
  author={Roessler, Damian},
  title={On $\psi $-invariant subvarieties of semiabelian varieties and the Manin-Mumford conjecture},
  journal={J. Algebraic Geom.},
  volume={13},
  date={2004},
  number={4},
  pages={771--798},
  issn={1056-3911},
}

\bib{Poonen-Voloch-The-Brauer-Manin}{article}{
   author={Poonen, Bjorn},
   author={Voloch, Jos{\'e} Felipe},
   title={The Brauer-Manin obstruction for subvarieties of abelian varieties
   over function fields},
   journal={Ann. of Math. (2)},
   volume={171},
   date={2010},
   number={1},
   pages={511--532},
   issn={0003-486X},
   doi={10.4007/annals.2010.171.511},
}

\bib{Raynaud-Faisceaux}{book}{
  author={Raynaud, Michel},
  title={Faisceaux amples sur les sch\'emas en groupes et les espaces homog\`enes},
  language={French},
  series={Lecture Notes in Mathematics, Vol. 119},
  publisher={Springer-Verlag},
  place={Berlin},
  date={1970},
  pages={ii+218},
}

\bib{Raynaud-Sous-var}{article}{
  author={Raynaud, Michel},
  title={Sous-vari\'et\'es d'une vari\'et\'e ab\'elienne et points de torsion},
  conference={ title={Arithmetic and geometry, Vol. I}, },
  book={ series={Progr. Math.}, volume={35}, publisher={Birkh\"auser Boston}, place={Boston, MA}, },
  date={1983},
  pages={327--352},
}

\bib{Rossler-On-the-Manin}{article}{
  author={R\"ossler, Damian},
  title={On the Manin-Mumford and Mordell-Lang conjectures in positive characteristic},
  status={ArXiv 1104.4311. Submitted}
}


\bib{Tate-p-divisible}{article}{
  author={Tate, John},
  title={$p$-divisible groups.},
  conference={ title={Proc. Conf. Local Fields}, address={Driebergen}, date={1966}, },
  book={ publisher={Springer}, place={Berlin}, },
  date={1967},
  pages={158--183},
}

\bib{Tate-Voloch-Linear}{article}{
  author={Tate, John},
  author={Voloch, Jos{\'e} Felipe},
  title={Linear forms in $p$-adic roots of unity},
  journal={Internat. Math. Res. Notices},
  date={1996},
  number={12},
  pages={589--601},
  issn={1073-7928},
}


\bib{Voloch-Diophantine}{article}{
  author={Voloch, Jos{\'e} Felipe},
  title={Diophantine approximation on abelian varieties in characteristic $p$},
  journal={Amer. J. Math.},
  volume={117},
  date={1995},
  number={4},
  pages={1089--1095},
  issn={0002-9327},
  doi={10.2307/2374961},
}



\end{biblist}
\end{bibdiv}
\end{document}